\documentclass[11pt]{article}
\usepackage{amssymb}
\usepackage{amsfonts}
\usepackage{graphicx}
\usepackage{amsmath}

\input tcilatex
\begin{document}

\begin{center}
\textbf{LACUNARY STATISTICAL CONVERGENCE OF ORDER }$\beta $

\textbf{IN DIFFERENCE SEQUENCES OF FUZZY NUMBERS}

\textit{\bigskip }

\textit{H\i fs\i\ ALTINOK\qquad Damla YA\u{G}DIRAN}

\textit{Department of Mathematics, F\i rat University, 23119, Elaz\i \u{g},
TURKEY.}

\textit{E-mail:} \textit{hifsialtinok@gmail.com\quad dyagdiran@hotmail.com}

\bigskip
\end{center}

\begin{quote}
\textbf{Abstract.} In this paper, we define the spaces $N_{\theta }^{\beta
}\left( p,F,\Delta ^{m}\right) ,$ $S_{\theta }^{\beta }\left( F,\Delta
^{m}\right) ,$ $w_{p}^{\beta }\left( F,\Delta ^{m}\right) $ for sequences of
fuzzy numbers using generalized difference operator $\Delta ^{m}$ and a
lacunary sequence $\theta $ and give some relations between them, where $%
\beta \in \left( 0,1\right] $ and $p>0$. Furthermore, in the last section of
paper, some inclusion theorems are presented related to the spaces $%
S_{\theta }^{\beta }\left( F,\Delta ^{m}\right) $ and $w_{p}^{\beta }\left(
\theta ,f,F,\Delta ^{m}\right) $ according to modulus function $f$.
\end{quote}

\bigskip



\bigskip

\noindent \textbf{1 Introduction}

\QTP{Body Math}
$\medskip $

The idea of statistical convergence was given by Zygmund \cite{zygmund} in
the first edition of his monograph published in Warsaw in 1935. The concept
of statistical convergence was introduced by Steinhaus \cite{steinhaus} and
Fast \cite{fast} and later reintroduced by Schoenberg \cite{schonberg}
independently. Over the years and under different names statistical
convergence has been discussed in the theory of Fourier analysis, ergodic
theory, number theory, measure theory, trigonometric series, turnpike theory
and Banach spaces. Later on it was further investigated from the sequence
space point of view and linked with summability theory by Connor \cite%
{connor}, Et \cite{et1}, Fridy \cite{fridy}, Mursaleen \cite{musaleen} and
many others. In recent years, generalizations of statistical convergence
have appeared in the study of strong integral summability and the structure
of ideals of bounded continuous functions on locally compact spaces.
Statistical convergence and its generalizations are also connected with
subsets of the Stone-\v{C}ech compactification of the natural numbers.
Moreover, statistical convergence is closely related to the concept of
convergence in probability.

\medskip

Now, we recall some basic concepts which we shall use throughout the paper.

By a lacunary sequence we mean an increasing integer sequence $\theta
=\left( k_{r}\right) $ such that $k_{0}=0$ and $h_{r}=\left(
k_{r}-k_{r-1}\right) \rightarrow \infty $ as $r\rightarrow \infty .$
Throught this paper the intervals determined by $\theta $ will be denoted by 
$I_{r}=\left( k_{r-1},k_{r}\right] $ and the ratio $\frac{k_{r}}{k_{r-1}}$
will be abbreviated by $q_{r}.$ Recently lacunary sequence spaces have been
studied in (\cite{fredman},\cite{fridy-orhan},\cite{altinok1})

The order of statistical convergence of a sequence of numbers was given by
Gadjiev and Orhan \cite{gadjiev} and later on statistical convergence of
order $\alpha $ and strong $p-$Ces\`{a}ro summability of order $\alpha $\
was studied by \c{C}olak \cite{colak}. Lacunary statistical convergence of
order $\alpha $ in real number sequences was defined by \c{S}eng\"{u}l and
Et \cite{sengul}.

The difference spaces $\ell _{\infty }\left( \Delta \right) $, $c\left(
\Delta \right) $ and $c_{0}\left( \Delta \right) $, consisting of all real
valued sequences $x=\left( x_{k}\right) $ such that $\Delta x=\Delta
^{1}x=\left( x_{k}-x_{k+1}\right) $ in the sequence spaces $\ell _{\infty }$%
, $c$ and $c_{0}$, were defined by K\i zmaz \cite{Kizmaz}. The idea of
difference sequences was generalized by Et and \c{C}olak \cite{Et and Colak}%
, Altinok \cite{altinok3}, \c{C}olak \textit{et al. }\cite{colak2}, Tripathy
and Baruah \cite{tripathy1} and many others.

\bigskip

A fuzzy number is a mapping $X:\mathbb{R}\longrightarrow \lbrack 0,1]$ which
satisfies the following properties:

$i)$ $X$ is normal, i.e., there exists an $t_{0}\in \mathbb{R}$ such that $%
X(t_{0})=1;$

$ii)$ $X$ is fuzzy convex, i.e., for all $\lambda \in \lbrack 0,1]$ and all $%
s,t\in \mathbb{R}$ and $0\leq \lambda \leq 1,X(\lambda s+(1-\lambda )t)\geq
\min [X(s),X(t)];$

$iii)$ $X$ is upper semicontinuous on real numbers set;

$iv)$ The set $\overline{\{t\in \mathbb{R}:X(t)>0\}},$ denoted by $[X]^{0},$
is compact, where $\overline{\{t\in \mathbb{R}:X(t)>0\}}$ is closure of $%
\{t\in \mathbb{R}:X(t)>0\}$ with usual topology of $\mathbb{R}$.

\noindent $L(\mathbb{R})$ denotes the set of all fuzzy numbers on real
numbers set $\mathbb{R}$ and is said to be a fuzzy number space.

For $\alpha \in \left( 0,1\right] $, the $\alpha $-level set $[X]^{\alpha }$
of fuzzy number $X$ is defined by$.$%
\begin{equation*}
\lbrack X]^{\alpha }=\left\{ 
\begin{array}{cc}
\{t\in \mathbb{R}:X(t)\geq \alpha \}, & \text{for }\alpha \in \left( 0,1%
\right] \\ 
\overline{\{t\in \mathbb{R}:X(t)>\alpha \}}, & \text{for }\alpha =0%
\end{array}%
\right.
\end{equation*}

The aritmetic operations for $\alpha $-level sets are defined as follows:

For $X,Y\in L(\mathbb{R}),$ now and what follows we will denote $\alpha -$%
level sets by $\left[ X\right] ^{\alpha }=\left[ u_{1}^{\alpha
},v_{1}^{\alpha }\right] ,$ $\left[ Y\right] ^{\alpha }=\left[ u_{2}^{\alpha
},v_{2}^{\alpha }\right] ,$ $\alpha \in \left[ 0,1\right] .$ Then we have%
\begin{align*}
\left[ X+Y\right] ^{\alpha }& =\left[ u_{1}^{\alpha }+u_{2}^{\alpha
},v_{1}^{\alpha }+v_{2}^{\alpha }\right] \\
\left[ X-Y\right] ^{\alpha }& =\left[ u_{1}^{\alpha }-v_{2}^{\alpha
},v_{1}^{\alpha }-u_{2}^{\alpha }\right] \\
\left[ X.Y\right] ^{\alpha }& =\left[ \min_{i,j\in \left[ 1,2\right]
}u_{i}^{\alpha }v_{j}^{\alpha },\max_{i,j\in \left[ 1,2\right]
}u_{i}^{\alpha }v_{j}^{\alpha }\right] .
\end{align*}

A sequence $X=(X_{k})$ of fuzzy numbers is a function $X$ from the set $%
\mathbb{N}$ of all natural numbers into $L(\mathbb{R})$ \cite{matloka}.

Let $\beta \in \left( 0,1\right] $ and $X=\left( X_{k}\right) $ be a
sequence of fuzzy numbers. Then the sequence $X=\left( X_{k}\right) $ of
fuzzy numbers is said to be statistically convergent of order $\beta ,$ to
fuzzy number $X_{0}$ if for every $\varepsilon >0,$%
\begin{equation*}
\lim_{n\rightarrow \infty }\frac{1}{n^{\beta }}\left\vert \left\{ k\leq
n:d\left( X_{k},X_{0}\right) \geq \varepsilon \right\} \right\vert =0,
\end{equation*}%
where the vertical bars indicate the number of elements in the enclosed set.
In this case we write $S^{\beta }\left( F\right) -\lim X_{k}=X_{0}.$ We
denote the set of all statistically convergent sequences of order $\beta $
by $S^{\beta }\left( F\right) $ \cite{altinok etal}.

Let $w\left( F\right) $ be the set of all sequences of fuzzy numbers. The
operator $\Delta ^{m}:w\left( F\right) \rightarrow w\left( F\right) $ is
defined by

\begin{eqnarray*}
\left( \Delta ^{0}X\right) _{k} &=&X_{k},\left( \Delta ^{1}X\right)
_{k}=\Delta ^{1}X_{k}=X_{k}-X_{k+1}, \\
\left( \Delta ^{m}X\right) _{k} &=&\Delta ^{1}\left( \Delta ^{m-1}X\right)
_{k},\left( m\geq 2\right) ,\text{ for all }k\in \mathbb{N}.
\end{eqnarray*}

\bigskip

\noindent \textbf{2 Main Results}

\bigskip

In this section, we define the spaces $S_{\theta }^{\beta }\left( F,\Delta
^{m}\right) $ and $N_{\theta }^{\beta }\left( p,F,\Delta ^{m}\right) $ and
examine some inclusion relations between them and space $w_{p}^{\beta
}\left( F,\Delta ^{m}\right) .$

\medskip

\noindent \textbf{Definition 2.1}\ Let $\theta =\left( k_{r}\right) $ be a
lacunary sequence, $X=\left( X_{k}\right) $ be a sequence of fuzzy numbers
and $\beta \in \left( 0,1\right] $ be given. The sequence $X=\left(
X_{k}\right) \in w\left( F\right) $ is said to be $S_{\theta }^{\beta
}\left( F,\Delta ^{m}\right) -$statistically convergent (or lacunary
statistically convergent sequence of order $\beta )$ if there is a fuzzy
number $X_{0}$ such that 
\begin{equation*}
\lim_{r\rightarrow \infty }\frac{1}{h_{r}^{\beta }}\left\vert \left\{ k\in
I_{r}:d\left( \Delta ^{m}X_{k},X_{0}\right) \geq \varepsilon \right\}
\right\vert =0,
\end{equation*}%
where $I_{r}=\left( k_{r-1},k_{r}\right] $ and $h_{r}^{\beta }$ denote the $%
\beta $th power $\left( h_{r}\right) ^{\beta }$ of $h_{r},$ that is $%
h^{\beta }=(h_{r}^{\beta })=(h_{1}^{\beta },h_{2}^{\beta },...,h_{r}^{\beta
},...).$ In this case we write $S_{\theta }^{\beta }\left( F,\Delta
^{m}\right) -\lim X_{k}=X_{0}$. The set of all $S_{\theta }^{\beta }\left(
F,\Delta ^{m}\right) -$statistically convergent sequences will be denoted by 
$S_{\theta }^{\beta }\left( F,\Delta ^{m}\right) $. For $\theta =\left(
2^{r}\right) ,$ we shall write $S^{\beta }\left( F,\Delta ^{m}\right) $
instead of $S_{\theta }^{\beta }\left( F,\Delta ^{m}\right) $ and we shall
write $S\left( F,\Delta ^{m}\right) $ instead of $S_{\theta }^{\beta }\left(
F,\Delta ^{m}\right) $ in the special case $\beta =1$ and $\theta =\left(
2^{r}\right) .$

\bigskip

\noindent \textbf{Lemma 2.2. }Let $\theta =\left( k_{r}\right) $ be a
lacunary sequence and $K\subset \mathbb{N}.$ Then $S_{\theta }^{\beta
}\left( K\right) \leq S_{\theta }^{\alpha }\left( K\right) $ for $0<\alpha
\leq \beta \leq 1$.

\textbf{Proof. }Let\textbf{\ }$0<\alpha \leq \beta \leq 1$. We can write

\begin{equation*}
\dfrac{1}{h_{r}^{\beta }}\left\vert \left\{ k_{r-1}<k\leq k_{r}:k\in
K\right\} \right\vert \leq \dfrac{1}{h_{r}^{\alpha }}\left\vert \left\{
k_{r-1}<k\leq k_{r}:k\in K\right\} \right\vert
\end{equation*}%
since $h_{r}^{\alpha }\leq h_{r}^{\beta }$. Hence $S_{\theta }^{\beta
}\left( K\right) \leq $ $S_{\theta }^{\alpha }\left( K\right) $.

\bigskip

Let $\beta \in \left( 0,1\right] .$ If a sequence of fuzzy numbers is
lacunary statistically convergent of order $\beta $, then it is lacunary
statistically convergent. Really, we obtain%
\begin{equation*}
\dfrac{1}{h_{r}^{\alpha }}\left\vert \left\{ k\in I_{r}:d\left( \Delta
^{m}X_{k},X_{0}\right) \geq \varepsilon \right\} \right\vert \geq \dfrac{1}{%
h_{r}}\left\vert \left\{ k\in I_{r}:d\left( \Delta ^{m}X_{k},X_{0}\right)
\geq \varepsilon \right\} \right\vert
\end{equation*}%
since $h_{r}^{\alpha }\leq h_{r}$ for $\beta \in \left( 0,1\right] .$ Taking
limit as $r\rightarrow \infty $ we get $X_{k}\rightarrow X_{0}\left(
S_{\theta }^{\alpha },\Delta ^{m}\right) $ and so $X_{k}\rightarrow
X_{0}\left( S_{\theta },\Delta ^{m}\right) .$

\medskip

Note that the lacunary statistical convergence of order $\beta $ is well
defined for $\beta \in \left( 0,1\right] $, but not well defined for $\beta
>1$. To show this, consider the sequence $X=\left( X_{k}\right) $ of fuzzy
numbers following:%
\begin{equation*}
X_{k}\left( t\right) =\left\{ 
\begin{array}{cc}
\left. 
\begin{array}{cc}
t, & \text{for }0\leq t\leq 1 \\ 
-t+2, & \text{for }1\leq t\leq 2 \\ 
0, & \text{otherwise}%
\end{array}%
\right\} :=X^{\prime }, & 
\begin{array}{c}
\text{if }k=2r \\ 
(r=1,2,3,...)%
\end{array}
\\ 
\left. 
\begin{array}{cc}
t-3, & \text{for }3\leq t\leq 4 \\ 
-t+5, & \text{for }4\leq t\leq 5 \\ 
0, & \text{otherwise}%
\end{array}%
\right\} :=X^{\prime \prime }, & \text{if }k\neq 2r\text{ }%
\end{array}%
\right.
\end{equation*}%
We can find the $\alpha -$level sets of sequences $\left( X_{k}\right) $ and 
$\left( \Delta ^{m}X_{k}\right) $ after some arithmetic operations as follows%
\begin{equation*}
\left[ X_{k}\right] ^{\alpha }=\left\{ 
\begin{array}{cc}
\left[ \alpha ,2-\alpha \right] , & \text{if }k=2r \\ 
\left[ 3+\alpha ,5-\alpha \right] , & \text{if }k\neq 2r\text{ }%
\end{array}%
\right. .
\end{equation*}%
and%
\begin{equation*}
\left[ \Delta ^{m}X_{k}\right] ^{\alpha }=\left\{ 
\begin{array}{cc}
\left[ 2^{m}\alpha -2^{m-1}.5,-2^{m}\alpha -2^{m-1}\right] :=X^{\prime }, & 
\text{if }k=2r \\ 
\left[ 2^{m}\alpha +2^{m-1},-2^{m}\alpha +2^{m-1}.5\right] :=X^{\prime
\prime }, & \text{if }k\neq 2r\text{ }%
\end{array}%
\right. ,\left( m=1,2,...\right)
\end{equation*}%
then we can write%
\begin{equation*}
\lim_{r\rightarrow \infty }\frac{1}{h_{r}^{\beta }}\left\vert \left\{ k\in
I_{r}:d\left( \Delta ^{m}X_{k},X^{\prime }\right) \geq \varepsilon \right\}
\right\vert \leq \lim_{r\rightarrow \infty }\frac{k_{r}-k_{r-1}}{%
2h_{r}^{\beta }}=\lim_{r\rightarrow \infty }\frac{h_{r}}{2h_{r}^{\beta }}=0
\end{equation*}%
and%
\begin{equation*}
\lim_{r\rightarrow \infty }\frac{1}{h_{r}^{\beta }}\left\vert \left\{ k\in
I_{r}:d\left( \Delta ^{m}X_{k},X^{\prime \prime }\right) \geq \varepsilon
\right\} \right\vert \leq \lim_{r\rightarrow \infty }\frac{k_{r}-k_{r-1}}{%
2h_{r}^{\beta }}=\lim_{r\rightarrow \infty }\frac{h_{r}}{2h_{r}^{\beta }}=0
\end{equation*}%
for $\beta >1.$ Hence sequence $\left( X_{k}\right) $ is $\Delta ^{m}-$%
lacunary statistically convergent of order $\beta ,$ both to $X^{\prime }$
and $X^{\prime \prime },$ i.e., $S_{\theta }^{\beta }\left( F,\Delta
^{m}\right) -\lim X_{k}$\ $=X^{\prime }$ and $S_{\theta }^{\beta }\left(
F,\Delta ^{m}\right) -\lim X_{k}$\ $=X^{\prime \prime }.$ But this is
impossible.

\medskip

\noindent \textbf{Theorem 2.3 }Let $0<\beta \leq 1$ and $X=\left(
X_{k}\right) ,$ $Y=\left( Y_{k}\right) $\ be sequences of fuzzy numbers$,$
then

$\left( i\right) $ If $S_{\theta }^{\beta }\left( F,\Delta ^{m}\right) -\lim
X_{k}=X_{0}$ and $c\in \mathbb{C},$ then $S_{\theta }^{\beta }\left(
F,\Delta ^{m}\right) -\lim \left( cx_{k}\right) =cx_{0},$

$\left( ii\right) $ If $S_{\theta }^{\beta }\left( F,\Delta ^{m}\right)
-\lim X_{k}=X_{0}$ and $S_{\theta }^{\beta }\left( F,\Delta ^{m}\right)
-\lim Y_{k}=Y_{0},$ then $S_{\theta }^{\beta }\left( F,\Delta ^{m}\right)
-\lim \left( X_{k}+Y_{k}\right) =X_{0}+Y_{0}.$

\textbf{Proof.} $(i)$ It is clear for the case $c=0.$ Suppose that $c\neq 0,$
then the proof of $\left( i\right) $ follows from%
\begin{equation*}
\frac{1}{h_{r}^{\beta }}\left\vert \left\{ k\in I_{r}:d\left( c\Delta
^{m}X_{k},cX_{0}\right) \geq \varepsilon \right\} \right\vert =\frac{1}{%
h_{r}^{\beta }}\left\vert \left\{ k\in I_{r}:d\left( \Delta
^{m}X_{k},X_{0}\right) \geq \frac{\varepsilon }{\left\vert c\right\vert }%
\right\} \right\vert
\end{equation*}%
and that of $\left( ii\right) $ follows from%
\begin{eqnarray*}
\frac{1}{h_{r}^{\beta }}\left\vert \left\{ k\in I_{r}:d\left( \Delta
^{m}\left( X_{k}+Y_{k}\right) ,\left( X_{0}+Y_{0}\right) \right) \geq
\varepsilon \right\} \right\vert &\leq &\frac{1}{h_{r}^{\beta }}\left\vert
\left\{ k\in I_{r}:d\left( \Delta ^{m}X_{k},X_{0}\right) \geq \frac{%
\varepsilon }{2}\right\} \right\vert \\
&&+\frac{1}{h_{r}^{\beta }}\left\vert \left\{ k\in I_{r}:d\left( \Delta
^{m}Y_{k},Y_{0}\right) \geq \frac{\varepsilon }{2}\right\} \right\vert .
\end{eqnarray*}

\noindent \textbf{Definition 2.4}\ Let $X=\left( X_{k}\right) $ be a
sequence of fuzzy numbers, $\theta =\left( k_{r}\right) $ be a lacunary
sequence, $\beta \in \left( 0,1\right] $ be any real number and let $p$ be a
positive real number. A sequence $X$ of fuzzy numbers is said to be strongly 
$N_{\theta }^{\beta }\left( p,F,\Delta ^{m}\right) -$summable (or strongly $%
N_{\theta }\left( p,F,\Delta ^{m}\right) -$summable of order $\beta )$ if
there is a fuzzy number $X_{0}$ such that 
\begin{equation*}
\underset{r\rightarrow \infty }{\lim }\frac{1}{h_{r}^{\beta }}%
\tsum\limits_{k\in I_{r}}d\left( \Delta ^{m}X_{k},X_{0}\right) ^{p}=0\text{.}
\end{equation*}%
In this case we write $N_{\theta }^{\beta }\left( p,F,\Delta ^{m}\right)
-\lim X_{k}=X_{0}$. The strong $N_{\theta }^{\beta }\left( p,F,\Delta
^{m}\right) -$summability reduces to the strong $N_{\theta }\left(
p,F,\Delta ^{m}\right) -$summability for $\beta =1$. The set of all strongly 
$N_{\theta }^{\beta }\left( p,F,\Delta ^{m}\right) -$summable sequences will
be denoted by $N_{\theta }^{\beta }\left( p,F,\Delta ^{m}\right) $.

If we take $\theta =\left( 2^{r}\right) $ in space $N_{\theta }^{\beta
}\left( p,F,\Delta ^{m}\right) ,$ then we obtain strongly $p-$Cesaro
summable sequences set of order $\beta $ following

\begin{equation*}
w_{p}^{\beta }\left( F,\Delta ^{m}\right) =\left\{ X=\left( X_{k}\right)
:\exists X_{0},\underset{n}{\lim }\dfrac{1}{n^{\beta }}\sum_{k=1}^{n}\left[
d\left( \Delta ^{m}X_{k},X_{0}\right) \right] ^{p}=0\right\}
\end{equation*}%
This set will be denoted by $w_{p,0}^{\beta }$ in case of $X_{0}=\bar{0}$ 
\cite{altinok3}.

\medskip

\noindent \textbf{Theorem 2.5}\ If $0<\beta <\gamma \leq 1$ then $S_{\theta
}^{\beta }\left( F,\Delta ^{m}\right) \subset S_{\theta }^{\gamma }\left(
F,\Delta ^{m}\right) $ and the inclusion is strict.

\medskip

\textbf{Proof. }Proof can be obtained from the inequality 
\begin{equation*}
\frac{1}{h_{r}^{\gamma }}\left\vert \left\{ k\in I_{r}:d\left( \Delta
^{m}X_{k},X_{0}\right) \geq \varepsilon \right\} \right\vert \leq \frac{1}{%
h_{r}^{\beta }}\left\vert \left\{ k\in I_{r}:d\left( \Delta
^{m}X_{k},X_{0}\right) \geq \varepsilon \right\} \right\vert .
\end{equation*}

\bigskip

Taking $\theta =\left( 2^{r}\right) $ we show the strictness of the
inclusion $S_{\theta }^{\beta }\left( F,\Delta ^{m}\right) \subset S_{\theta
}^{\gamma }\left( F,\Delta ^{m}\right) $ for a special case. For this,
consider the sequence $X=\left( X_{k}\right) $ of fuzzy numbers defined by%
\begin{equation*}
X_{k}\left( t\right) =\left\{ 
\begin{array}{cc}
\left. 
\begin{array}{cc}
2t-\left( 2k-1\right) , & \text{for }k-\frac{1}{2}\leq t\leq k \\ 
-2t+\left( 2k+1\right) , & \text{for }k\leq t\leq k+\frac{1}{2} \\ 
0, & \text{otherwise}%
\end{array}%
\right\} , & \text{if }k=n^{3} \\ 
\left. 
\begin{array}{cc}
2t-1, & \text{for }\frac{1}{2}\leq t\leq 1 \\ 
-2t+3, & \text{for }1\leq t\leq \frac{3}{2} \\ 
0, & \text{otherwise}%
\end{array}%
\right\} , & \text{if }k\neq n^{3}%
\end{array}%
\right.
\end{equation*}%
After calculating $\alpha -$level sets of sequences $\left( X_{k}\right) $
and $\left( \Delta ^{m}X_{k}\right) ,$ we obtain $S_{\theta }^{\gamma
}\left( F,\Delta ^{m}\right) -\lim X_{k}=\left[ 2^{m-1}\left( \alpha
-1\right) ,2^{m-1}\left( 1-\alpha \right) \right] ,$ i.e. $X\in S_{\theta
}^{\gamma }\left( F,\Delta ^{m}\right) $ for $\frac{1}{3}<\gamma \leq 1,$
but $X\notin S_{\theta }^{\beta }\left( F,\Delta ^{m}\right) $ for $0<\beta
\leq \frac{1}{3}.$ This implies that the inclusion $S_{\theta }^{\beta
}\left( F,\Delta ^{m}\right) \subset S_{\theta }^{\gamma }\left( F,\Delta
^{m}\right) $ is strict for $\beta ,\gamma \in \left( 0,1\right] $ such that 
$\beta \in \left( 0,\frac{1}{3}\right] $ and $\gamma \in \left( \frac{1}{3},1%
\right] .$

\medskip

\noindent \textbf{Corollary 2.6} If a sequence of fuzzy numbers is $%
S_{\theta }^{\beta }\left( F,\Delta ^{m}\right) -$statistically convergent
to fuzzy number $X_{0}$, then it is $S_{\theta }\left( F,\Delta ^{m}\right)
- $statistically convergent to fuzzy number $X_{0}.$

\medskip

\noindent \textbf{Theorem 2.7} Let $\beta $ and $\gamma $ be fixed real
numbers such that $0<\beta \leq \gamma \leq 1,$ $X=\left( X_{k}\right) $ be
a sequence of fuzzy numbers and $\theta =\left( k_{r}\right) $ be a lacunary
sequence. For $0<p<\infty ,$ $N_{\theta }^{\beta }\left( p,F,\Delta
^{m}\right) \subset S_{\theta }^{\gamma }\left( F,\Delta ^{m}\right) $ and
the inclusion is strict for some $\beta $'s and $\gamma $'s$.$

\medskip

\textbf{Proof.} For any sequence $X=\left( X_{k}\right) $ of fuzzy numbers
and $\varepsilon >0$, we have%
\begin{eqnarray*}
\tsum\limits_{k\in I_{r}}\left[ d\left( \Delta ^{m}X_{k},X_{0}\right) \right]
^{p} &=&\tsum\limits_{\substack{ k\in I_{r}  \\ d\left( \Delta
^{m}X_{k},X_{0}\right) \geq \varepsilon }}\left[ d\left( \Delta
^{m}X_{k},X_{0}\right) \right] ^{p}+\tsum\limits_{\substack{ k\in I_{r}  \\ %
d\left( \Delta ^{m}X_{k},X_{0}\right) <\varepsilon }}\left[ d\left( \Delta
^{m}X_{k},X_{0}\right) \right] ^{p} \\
&\geq &\tsum\limits_{\substack{ k\in I_{r}  \\ d\left( \Delta
^{m}X_{k},X_{0}\right) \geq \varepsilon }}\left[ d\left( \Delta
^{m}X_{k},X_{0}\right) \right] ^{p} \\
&\geq &\left\vert \left\{ k\in I_{r}:d\left( \Delta ^{m}X_{k},X_{0}\right)
\geq \varepsilon \right\} \right\vert \varepsilon ^{p}
\end{eqnarray*}%
and so that 
\begin{eqnarray*}
\frac{1}{h_{r}^{\beta }}\tsum\limits_{k\in I_{r}}\left[ d\left( \Delta
^{m}X_{k},X_{0}\right) \right] ^{p} &\geq &\frac{1}{h_{r}^{\beta }}%
\left\vert \left\{ k\in I_{r}:d\left( \Delta ^{m}X_{k},X_{0}\right) \geq
\varepsilon \right\} \right\vert \varepsilon ^{p} \\
&\geq &\frac{1}{h_{r}^{\gamma }}\left\vert \left\{ k\in I_{r}:d\left( \Delta
^{m}X_{k},X_{0}\right) \geq \varepsilon \right\} \right\vert \varepsilon
^{p}.
\end{eqnarray*}%
It follows that if $X=\left( X_{k}\right) $ is strongly $N_{\theta }^{\beta
}\left( p,F,\Delta ^{m}\right) -$summable to fuzzy number $X_{0}$, then it
is $S_{\theta }^{\gamma }\left( F,\Delta ^{m}\right) -$statistically
convergent to $X_{0}$.

\medskip

Taking $\beta =\gamma ,$ $\theta =\left( 2^{r}\right) $ and $p=1$ we show
the strictness of the inclusion $N_{\theta }^{\beta }\left( p,F,\Delta
^{m}\right) \subset S_{\theta }^{\gamma }\left( F,\Delta ^{m}\right) $ for a
special case. For this we can select a sequence of fuzzy numbers as follows:

\begin{equation*}
X_{k}\left( t\right) =\left\{ 
\begin{array}{cc}
\left. 
\begin{array}{cc}
\frac{t}{k}+1, & \text{for }-k\leq t\leq 0 \\ 
-\frac{t}{k}+1, & \text{for }0\leq t\leq k \\ 
0, & \text{otherwise}%
\end{array}%
\right\} , & 
\begin{array}{c}
\text{if }k=n^{2} \\ 
\left( n=1,2,...\right)%
\end{array}
\\ 
\left. 
\begin{array}{cc}
\frac{t}{2}-1, & \text{for }2\leq t\leq 4 \\ 
-\frac{t}{2}+3, & \text{for }4\leq t\leq 6 \\ 
0, & \text{otherwise}%
\end{array}%
\right\} :=X_{0} & \text{if }k\neq n^{2}%
\end{array}%
\right.
\end{equation*}%
After calculating $\alpha -$level sets of sequences $\left( X_{k}\right) $
and $\left( \Delta ^{m}X_{k}\right) ,$ we see that $\left( X_{k}\right) $ is 
$\Delta ^{m}-$statistically convergent of order $\gamma $ to $X_{0}$ for $%
\frac{1}{2}<\beta \leq 1,$ where $\left[ X_{0}\right] ^{\alpha }=\left[
2^{m+1}\left( \alpha -1\right) ,2^{m}\left( 1-\alpha \right) \right] ,$ but
it is not strongly $N_{\theta }^{\beta }\left( p,F,\Delta ^{m}\right) $%
-summable to $X_{0}.$ Therefore $N_{\theta }^{\beta }\left( p,F,\Delta
^{m}\right) \subset S_{\theta }^{\beta }\left( F,\Delta ^{m}\right) $ for $%
\frac{1}{2}<\beta \leq 1.$

\medskip

\noindent \textbf{Corollary 2.8 }Let $X=\left( X_{k}\right) $ be a sequence
of fuzzy numbers, $\beta \in \left( 0,1\right] $ and $p$ is a positive real
number. Then

$i)$ If $X\in N_{\theta }^{\beta }\left( p,F,\Delta ^{m}\right) ,$ then $%
X\in S_{\theta }^{\beta }\left( F,\Delta ^{m}\right) $ and the limits are
same.

$ii)$ If $X\in N_{\theta }^{\beta }\left( p,F,\Delta ^{m}\right) ,$ then $%
X\in S_{\theta }\left( F,\Delta ^{m}\right) $ and the limits are same.

\medskip

\noindent \textbf{Theorem 2.9 }Let $X=\left( X_{k}\right) $ be a sequence of
fuzzy numbers, $\beta \in \left( 0,1\right] $ and $\theta =\left(
k_{r}\right) $ be a lacunary sequence. If $\lim \inf_{r}q_{r}>1,$ then $%
S^{\beta }\left( F,\Delta ^{m}\right) \subset S_{\theta }^{\beta }\left(
F,\Delta ^{m}\right) .$

\textbf{Proof. }Suppose that $\lim \inf_{r}q_{r}>1,$ then there exists a
number $\delta >0$ such that $q_{r}\geq 1+\delta $ for sufficiently large $%
r\ $and so%
\begin{equation*}
\frac{h_{r}}{k_{r}}\geq \frac{\delta }{1+\delta }\Longrightarrow \left( 
\frac{h_{r}}{k_{r}}\right) ^{\beta }\geq \left( \frac{\delta }{1+\delta }%
\right) ^{\beta }\Longrightarrow \frac{1}{k_{r}^{\beta }}\geq \frac{\delta
^{\beta }}{\left( 1+\delta \right) ^{\beta }}\frac{1}{h_{r}^{\beta }}.
\end{equation*}%
If $X_{k}\rightarrow X_{0}\left( S^{\beta }\left( F,\Delta ^{m}\right)
\right) ,$ then we prove the suff\i ciency from following inequalities for
every $\varepsilon >0$ and for sufficiently large $r$ 
\begin{eqnarray*}
\frac{1}{k_{r}^{\beta }}\left\vert \left\{ k\leq k_{r}:d\left( \Delta
^{m}X_{k},X_{0}\right) \geq \varepsilon \right\} \right\vert &\geq &\frac{1}{%
k_{r}^{\beta }}\left\vert \left\{ k\in I_{r}:d\left( \Delta
^{m}X_{k},X_{0}\right) \geq \varepsilon \right\} \right\vert \\
&\geq &\frac{\delta ^{\beta }}{\left( 1+\delta \right) ^{\beta }}\frac{1}{%
h_{r}^{\beta }}\left\vert \left\{ k\in I_{r}:d\left( \Delta
^{m}X_{k},X_{0}\right) \geq \varepsilon \right\} \right\vert ;
\end{eqnarray*}

\noindent \textbf{Theorem 2.10} If $\lim_{r\rightarrow \infty }\inf \frac{%
h_{r}^{\beta }}{k_{r}}>0,$ then $S\left( F,\Delta ^{m}\right) \subset
S_{\theta }^{\beta }\left( F,\Delta ^{m}\right) .$

\textbf{Proof. }For a given $\varepsilon >0,$ we can write%
\begin{equation*}
\left\{ k\leq k_{r}:d\left( \Delta ^{m}X_{k},X_{0}\right) \geq \varepsilon
\right\} \supset \left\{ k\in I_{r}:d\left( \Delta ^{m}X_{k},X_{0}\right)
\geq \varepsilon \right\} .
\end{equation*}%
From here we obtain%
\begin{eqnarray*}
\frac{1}{k_{r}}\left\vert \left\{ k\leq k_{r}:d\left( \Delta
^{m}X_{k},X_{0}\right) \geq \varepsilon \right\} \right\vert &\geq &\frac{1}{%
k_{r}}\left\vert \left\{ k\in I_{r}:d\left( \Delta ^{m}X_{k},X_{0}\right)
\geq \varepsilon \right\} \right\vert \\
&=&\frac{h_{r}^{\beta }}{k_{r}}\frac{1}{h_{r}^{\beta }}\left\vert \left\{
k\in I_{r}:d\left( \Delta ^{m}X_{k},X_{0}\right) \geq \varepsilon \right\}
\right\vert .
\end{eqnarray*}%
Taking limit according to $r$ and using $\lim_{r\rightarrow \infty }\inf 
\frac{h_{r}^{\beta }}{k_{r}}>0,$ we get%
\begin{equation*}
X_{k}\rightarrow X_{0}\left( S,F,\Delta ^{m}\right) \Longrightarrow
X_{k}\rightarrow X_{0}\left( S_{\theta }^{\beta },F,\Delta ^{m}\right) .
\end{equation*}

\noindent \textbf{Theorem 2.11 }Let $X=\left( X_{k}\right) $ be a sequence
of fuzzy numbers, $0<\beta \leq 1,$and $\theta =\left( k_{r}\right) $ be a
lacunary sequence. If $\lim \sup_{r}q_{r}<\infty ,$ then $S_{\theta }^{\beta
}\left( F,\Delta ^{m}\right) \subset S\left( F,\Delta ^{m}\right) .$

\textbf{Proof. }Omitted.

\medskip

\noindent \textbf{Theorem 2.12}\ Let $0<\beta <\gamma \leq 1,$ $X=\left(
X_{k}\right) $ be a sequence of fuzzy numbers and $p>0$, then there exists
inclusion $N_{\theta }^{\beta }\left( p,F,\Delta ^{m}\right) \subset
N_{\theta }^{\gamma }\left( p,F,\Delta ^{m}\right) $ and this is strict.

\textbf{Proof. }Let $X=\left( X_{k}\right) \in N_{\theta }^{\beta }\left(
p,F,\Delta ^{m}\right) $. Then given $\beta $ and $\gamma $ such that $%
0<\beta <\gamma \leq 1$ and a $p>0,$ we can write inequality 
\begin{equation*}
\frac{1}{h_{r}^{\gamma }}\tsum\limits_{k\in I_{r}}\left[ d\left( \Delta
^{m}X_{k},X_{0}\right) \right] ^{p}\leq \frac{1}{h_{r}^{\beta }}%
\tsum\limits_{k\in I_{r}}\left[ d\left( \Delta ^{m}X_{k},X_{0}\right) \right]
^{p}
\end{equation*}%
and this gives that $N_{\theta }^{\beta }\left( p,F,\Delta ^{m}\right)
\subset N_{\theta }^{\gamma }\left( p,F,\Delta ^{m}\right) .$

\bigskip

To show that the inclusion is strict take lacunary sequence $\theta =\left(
2^{r}\right) $ and $p=1.$ Consider the sequence $X=\left( X_{k}\right) $ of
fuzzy numbers as follows 
\begin{equation*}
X_{k}\left( t\right) =\left\{ 
\begin{array}{cc}
\left. 
\begin{array}{cc}
t-2, & \text{for }2\leq t\leq 3 \\ 
-t+4, & \text{for }3\leq t\leq 4 \\ 
0, & \text{otherwise}%
\end{array}%
\right\} , & 
\begin{array}{c}
\text{if }k=n^{3} \\ 
\left( n=1,2,...\right)%
\end{array}
\\ 
\left. 
\begin{array}{cc}
\frac{t}{3}-\frac{5}{3}, & \text{for }5\leq t\leq 8 \\ 
-\frac{t}{3}+\frac{11}{3}, & \text{for }8\leq t\leq 11 \\ 
0, & \text{otherwise}%
\end{array}%
\right\} :=X_{0} & \text{if }k\neq n^{3}%
\end{array}%
\right.
\end{equation*}%
After calculating $\alpha -$level sets of sequences $\left( X_{k}\right) $
and $\left( \Delta ^{m}X_{k}\right) ,$ we obtain $X\in N_{\theta }^{\gamma
}\left( p,F,\Delta ^{m}\right) $ since $\left( X_{k}\right) $ is strongly $%
N_{\theta }^{\beta }\left( p,F,\Delta ^{m}\right) $-summable to $X_{0}$ for $%
\frac{1}{3}<\gamma \leq 1,$ where $\left[ X_{0}\right] ^{\alpha }=\left[
2^{m}\left( 3\alpha -3\right) ,2^{m}\left( 3-3\alpha \right) \right] ,$ but $%
X\notin N_{\theta }^{\beta }\left( p,F,\Delta ^{m}\right) $ for $0<\beta
\leq \frac{1}{3},$ i.e. it is not strongly $N_{\theta }^{\beta }\left(
p,F,\Delta ^{m}\right) $-summable to fuzzy number $X_{0}.$ So, this shows
that inclusion is strict. We show the limits of sequence $\Delta ^{m}X_{k}$
in Figure 1.%
\begin{equation*}
\FRAME{itbpF}{5.376in}{1.6835in}{0in}{}{}{Figure}{\special{language
"Scientific Word";type "GRAPHIC";maintain-aspect-ratio TRUE;display
"USEDEF";valid_file "T";width 5.376in;height 1.6835in;depth
0in;original-width 5.4496in;original-height 1.687in;cropleft "0";croptop
"1";cropright "1";cropbottom "0";tempfilename
'NUPD7N00.wmf';tempfile-properties "XPR";}}
\end{equation*}

\medskip

\noindent \textbf{Corollary 2.13 }Let $0<\beta \leq \gamma \leq 1,$ $%
X=\left( X_{k}\right) $ be a sequence of fuzzy numbers and $p>0$. Then

$(i)$ If $\beta =\gamma ,$ then $N_{\theta }^{\beta }\left( p,F,\Delta
^{m}\right) =N_{\theta }^{\gamma }\left( p,F,\Delta ^{m}\right) ,$

$(ii)$ For each $\beta \in (0,1]$ and $0<p<\infty ,$ $N_{\theta }^{\beta
}\left( p,F,\Delta ^{m}\right) \subseteq N_{\theta }\left( p,F,\Delta
^{m}\right) .$

\bigskip

\noindent \textbf{Theorem 2.14. }Let $X=\left( X_{k}\right) $ be a sequence
of fuzzy numbers and $\theta =\left( k_{r}\right) $ be a lacunary sequence.
If $\lim \inf_{r}q_{r}>1,$ then $w_{p}^{\beta }\left( F,\Delta ^{m}\right)
\subseteq N_{\theta }^{\beta }\left( p,F,\Delta ^{m}\right) $ for $0<\beta
\leq 1$ and $0<p<\infty .$

\textbf{Proof. }If $\lim \inf_{r}q_{r}>1,$ then there exists a number $%
\delta >0$ such that $q_{r}\geq 1+\delta $ for all $r\geq 1.$ For $0<\beta
\leq 1$ and $X\in w_{p,0}^{\beta }$ we can write

\begin{eqnarray*}
\tau _{r}^{\beta } &=&\dfrac{1}{h_{r}^{\beta }}\sum_{i=1}^{k_{r}}d\left(
\Delta ^{m}X_{i},\bar{0}\right) ^{p}-\dfrac{1}{h_{r}^{\beta }}%
\sum_{i=1}^{k_{r-1}}d\left( \Delta ^{m}X_{i},\bar{0}\right) ^{p} \\
&=&\dfrac{k_{r}^{\beta }}{h_{r}^{\beta }}\left( \dfrac{1}{k_{r}^{\beta }}%
\sum_{i=1}^{k_{r}}d\left( \Delta ^{m}X_{i},\bar{0}\right) ^{p}\right) -%
\dfrac{k_{r-1}^{\beta }}{h_{r}^{\beta }}\left( \dfrac{1}{k_{r-1}^{\beta }}%
\sum_{i=1}^{k_{r-1}}d\left( \Delta ^{m}X_{i},\bar{0}\right) ^{p}\right)
\end{eqnarray*}%
and so

\begin{equation*}
\dfrac{k_{r}^{\beta }}{h_{r}^{\beta }}\leq \dfrac{\left( 1+\delta \right)
^{\beta }}{\delta ^{\beta }}\text{ and }\dfrac{k_{r-1}^{\beta }}{%
h_{r}^{\beta }}\leq \dfrac{1}{\delta ^{\beta }}
\end{equation*}%
for $h_{r}=k_{r}-k_{r-1}$. It can be seen that both of terms$\dfrac{1}{%
k_{r}^{\beta }}\sum_{i=1}^{k_{r}}d\left( \Delta ^{m}X_{i},\bar{0}\right)
^{p} $ and $\dfrac{1}{k_{r-1}^{\beta }}\sum_{i=1}^{k_{r-1}}d\left( \Delta
^{m}X_{i},\bar{0}\right) ^{p}$ converge to $0$ and so $\tau _{r}^{\beta }$
does. Hence $X\in N_{\theta ,0}^{\beta }\left( p,F,\Delta ^{m}\right) $,
that is $w_{p}^{\beta }\left( F,\Delta ^{m}\right) \subseteq N_{\theta
}^{\beta }\left( p,F,\Delta ^{m}\right) $.

\bigskip

\noindent \textbf{Theorem 2.15.} If $\lim \sup_{r}\dfrac{k_{r}}{%
k_{r-1}^{\beta }}<\infty $, then $N_{\theta }\left( p,F,\Delta ^{m}\right)
\subseteq w_{p}^{\beta }\left( F,\Delta ^{m}\right) $ for $0<\beta \leq 1.$

\textbf{Proof. }If $\lim \sup_{r}\dfrac{k_{r}}{k_{r-1}^{\beta }}<\infty ,$
then there exists a number $M>0$ such that $\dfrac{k_{r}}{k_{r-1}^{\beta }}%
<M $ for all $r\geq 1.$ Let $X\in N_{\theta ,0}\left( p,F,\Delta ^{m}\right) 
$ and $\varepsilon >0.$ Then we can find numbers $R>0$ and $K>0$ such that $%
\sup_{i\geq R}\tau _{i}<\varepsilon $ and $\tau _{i}<K$ for $i=1,2,3,...$ .
Therefore we can write%
\begin{eqnarray*}
&&\dfrac{1}{t^{\beta }}\sum_{i=1}^{t}d\left( \Delta ^{m}X_{i},\bar{0}\right)
^{p}%
\begin{array}{c}
\leq%
\end{array}%
\dfrac{1}{k_{r-1}^{\beta }}\sum_{i=1}^{k_{r-1}}d\left( \Delta ^{m}X_{i},\bar{%
0}\right) ^{p} \\
&&\qquad \qquad 
\begin{array}{c}
=%
\end{array}%
\dfrac{1}{k_{r-1}^{\beta }}\left( \sum_{I_{1}}d\left( \Delta ^{m}X_{i},\bar{0%
}\right) ^{p}+\sum_{I_{2}}d\left( \Delta ^{m}X_{i},\bar{0}\right)
^{p}+...+\sum_{I_{r}}d\left( \Delta ^{m}X_{i},\bar{0}\right) ^{p}\right)
\end{eqnarray*}
\begin{eqnarray*}
&&\qquad \qquad 
\begin{array}{c}
=%
\end{array}%
\dfrac{k_{1}}{k_{r-1}^{\beta }}\tau _{1}+\dfrac{k_{2}-k_{1}}{k_{r-1}^{\beta }%
}\tau _{2}+...+\dfrac{k_{R}-k_{R-1}}{k_{r-1}^{\beta }}\tau _{R}+\dfrac{%
k_{R+1}-k_{R}}{k_{r-1}^{\beta }}\tau _{R+1}+...+\dfrac{kr-k_{r-1}}{%
k_{r-1}^{\beta }}\tau _{r} \\
&&\qquad \qquad 
\begin{array}{c}
\leq%
\end{array}%
\left( \sup_{i\geq 1}\tau _{i}\right) \dfrac{k_{R}}{k_{r-1}^{\beta }}+\left(
\sup_{i\geq R}\tau _{i}\right) \dfrac{k_{r}-k_{R}}{k_{r-1}^{\beta }}%
\begin{array}{c}
<%
\end{array}%
K\dfrac{k_{R}}{k_{r-1}^{\beta }}+\epsilon M
\end{eqnarray*}%
where $k_{r-1}<t\leq k_{r}$ and $r>R.$ We get $\dfrac{1}{t^{\beta }}%
\sum_{i=1}^{t}d\left( \Delta ^{m}X_{i},\bar{0}\right) ^{p}\rightarrow 0$
since $t\rightarrow \infty $ for $k_{r-1}\rightarrow \infty ,$ hence $X\in
w_{p,0}^{\beta }\left( F,\Delta ^{m}\right) .$

\bigskip

\noindent \textbf{Theorem 2.16. }Let $X=\left( X_{k}\right) $ be a sequence
of fuzzy numbers, $\theta =\left( k_{r}\right) $ be a lacunary sequence and $%
\beta \in \left( 0,1\right] $. If $X\in w^{\beta }\left( F,\Delta
^{m}\right) \cap N_{\theta }^{\beta }\left( p,F,\Delta ^{m}\right) $ and $%
\lim \sup_{r}\dfrac{k_{r}}{k_{r-1}^{\beta }}<\infty ,$ then the limits of $%
\left( X_{k}\right) $ are same.

\textbf{Proof. }Suppose that $w_{p}^{\beta }\left( F,\Delta ^{m}\right)
-\lim X_{k}=X^{\prime },$ $N_{\theta ,p}^{\beta }\left( p,F,\Delta
^{m}\right) -\lim X_{k}=X^{\prime \prime }$ and $X^{\prime }\neq X^{\prime
\prime }$. We get $N_{\theta ,0}\left( p,F,\Delta ^{m}\right) \subseteq
w_{p,0}^{\beta }\left( F,\Delta ^{m}\right) $ since $\lim \sup_{r}\dfrac{%
k_{r}}{k_{r-1}^{\beta }}<\infty $ from Theorem 2.15$.$ Because of $d\left(
\Delta ^{m}X_{i},X^{\prime \prime }\right) \in N_{\theta ,0}\left(
p,F,\Delta ^{m}\right) $, so $d\left( \Delta ^{m}X_{i},X^{\prime \prime
}\right) \in w_{p,0}^{\beta }\left( F,\Delta ^{m}\right) $. Therefore for $%
p=1$

\begin{equation*}
\dfrac{1}{t^{\beta }}\sum_{i=1}^{t}d\left( \Delta ^{m}X_{i},X^{\prime \prime
}\right) \rightarrow 0
\end{equation*}%
In the following inequality, both of terms at left converges to $0,$ but
this is a contradiction and so $X^{\prime }=X^{\prime \prime }.$

\begin{equation*}
\dfrac{1}{t^{\beta }}\sum_{i=1}^{t}d\left( \Delta ^{m}X_{i},X^{\prime \prime
}\right) +\dfrac{1}{t^{\beta }}\sum_{i=1}^{t}d\left( \Delta
^{m}X_{i},X^{\prime }\right) \geq \dfrac{1}{t^{\beta }}d\left( X^{\prime
},X^{\prime \prime }\right) >0
\end{equation*}

\bigskip

\noindent \textbf{3. Results Related to Modulus Function}

\bigskip

We give relation between sets $S_{\theta }^{\beta }\left( F,\Delta
^{m}\right) $ and $w_{p}^{\beta }\left[ \theta ,f,F,\Delta ^{m}\right] $
according to modulus function $f$.

\bigskip

We first quote the definition of a modulus function:

\bigskip

The concept of modulus function was formally introduced by Nakano \cite%
{nakano}. A mapping $f:[0,\infty )$ $\rightarrow $ $[0,\infty )$ is said to
be a modulus if

$i)$ $f\left( x\right) =0$ iff $x=0$,

$ii)$ $f\left( x+y\right) \leq f\left( x\right) +f\left( y\right) $ for $%
x,y\geq 0,$

$iii)$ $f$ is increasing,

$iv)$ $f$ is right-continuous at $0.$

The continuity of $f$ everywhere on $[0,\infty )$ follows from above
definition. A modulus function can be bounded or unbounded. For example $%
f\left( x\right) =x^{p},\left( 0<p\leq 1\right) $ is bounded and $f\left(
x\right) =\frac{x}{1+x}$ is bounded. The concept of modulus function for
sequences of fuzzy numbers was first investigated by Sarma \cite{Sarma}.
Later on it was explotied by Talo and Ba\c{s}ar \cite{talo2}{\footnotesize .}

\bigskip

\noindent \textbf{Definition 3.1. }Let $f$ be a modulus function, $X=\left(
X_{k}\right) $ be a sequence of fuzzy numbers and $p=\left( p_{k}\right) $
be a sequence of positive real numbers and $\beta \in \left( 0,1\right] $.
We define space $w_{p}^{\beta }\left( \theta ,f,F,\Delta ^{m}\right) $ as
follows

\begin{equation*}
w_{p}^{\beta }\left( \theta ,f,F,\Delta ^{m}\right) =\left\{ x=\left(
x_{k}\right) :\lim_{r\rightarrow \infty }\dfrac{1}{h_{r}^{\beta }}\sum_{k\in
I_{r}}\left[ d\left( \Delta ^{m}X_{k},X_{0}\right) \right] ^{p_{k}}=0,\text{
for some }X_{0}\right\}
\end{equation*}%
We shall write $N_{\theta }^{\beta }\left( p,F,\Delta ^{m}\right) $ instead
of $w_{p}^{\beta }\left( \theta ,f,F,\Delta ^{m}\right) $ in the special
cases $p_{k}=1$ and $f\left( x\right) =x$ for all $k\in 
\mathbb{N}
$.

We consider that $p=\left( p_{k}\right) $ is bounded and $%
0<h=\inf_{k}p_{k}\leq p_{k}\leq \sup_{k}p_{k}=H<\infty $ in the following
theorems.

\bigskip

\noindent \textbf{Theorem 3.2. }Let $\beta ,\gamma \in \left( 0,1\right] ,$ $%
\beta \leq \gamma ,$ $f$ be a modulus function, $X=\left( X_{k}\right) $ be
a sequence of fuzzy numbers and $\theta =\left( k_{r}\right) $ be a lacunary
sequence. Then, $w_{p}^{\beta }\left( \theta ,f,F,\Delta ^{m}\right) \subset
S_{\theta }^{\gamma }\left( F,\Delta ^{m}\right) .$

\textbf{Proof. }Take\textbf{\ }$x\in w_{p}^{\beta }\left( \theta ,f,F,\Delta
^{m}\right) $, $\varepsilon >0$ and let $\sum\limits_{1}$ and $%
\sum\limits_{2}$ be summations on $k\in I_{r},$ $d\left( \Delta
^{m}X_{k},X_{0}\right) \geq \varepsilon $ and $k\in I_{r},$ $d\left( \Delta
^{m}X_{k},X_{0}\right) <\varepsilon ,$ respectively. Since $h_{r}^{\beta
}<h_{r}^{\gamma }$ for each $r,$ we get

\begin{eqnarray*}
\dfrac{1}{h_{r}^{\beta }}\sum_{k\in I_{r}}\left[ f\left( d\left( \Delta
^{m}X_{k},X_{0}\right) \right) \right] ^{p_{k}} &=&\dfrac{1}{h_{r}^{\beta }}%
\left[ \sum_{1}\left[ f\left( d\left( \Delta ^{m}X_{k},X_{0}\right) \right) %
\right] ^{p_{k}}+\sum_{2}\left[ f\left( d\left( \Delta
^{m}X_{k},X_{0}\right) \right) \right] ^{p_{k}}\right] \\
&\geq &\dfrac{1}{h_{r}^{\gamma }}\left[ \sum_{1}\left[ f\left( d\left(
\Delta ^{m}X_{k},X_{0}\right) \right) \right] ^{p_{k}}+\sum_{2}\left[
f\left( d\left( \Delta ^{m}X_{k},X_{0}\right) \right) \right] ^{p_{k}}\right]
\\
&\geq &\dfrac{1}{h_{r}^{\gamma }}\sum_{1}\left[ f\left( \varepsilon \right) %
\right] ^{p_{k}} \\
&\geq &\dfrac{1}{h_{r}^{\gamma }}\sum_{1}\min \left( \left[ f\left(
\varepsilon \right) \right] ^{h},\left[ f\left( \varepsilon \right) \right]
^{H}\right) \\
&\geq &\dfrac{1}{h_{r}^{\gamma }}\left\vert \left\{ k\in I_{r}:d\left(
\Delta ^{m}X_{k},X_{0}\right) \geq \varepsilon \right\} \right\vert \min
\left( \left[ f\left( \varepsilon \right) \right] ^{h},\left[ f\left(
\varepsilon \right) \right] ^{H}\right)
\end{eqnarray*}%
Left side of above inequality converges to $0$ as \ $r\rightarrow \infty $
because of $X\in w_{p}^{\beta }\left( \theta ,f,F,\Delta ^{m}\right) $, so
right one converges to zero and hence $X\in S_{\theta }^{\gamma }\left(
F,\Delta ^{m}\right) $.

\bigskip

\noindent \textbf{Theorem 3.3. }Let $f$ be a bounded modulus function, $%
X=\left( X_{k}\right) $ be a sequence of fuzzy numbers and $\theta =\left(
k_{r}\right) $ be a lacunary sequence. If $\lim\limits_{r\rightarrow \infty }%
\dfrac{h_{r}}{h_{r}^{\beta }}=1,$ then $S_{\theta }^{\beta }\left( F,\Delta
^{m}\right) \subset w_{p}^{\beta }\left( \theta ,f,F,\Delta ^{m}\right) $.

\bigskip

\textbf{Proof. }Suppose that $f$ is bounded and $X\in S_{\theta }^{\beta
}\left( F,\Delta ^{m}\right) $. Given $\varepsilon >0.$\ There is an integer 
$K$ such that $\left( x\right) \leq K$ since $f$ is bounded. Now we can write

\begin{eqnarray*}
\dfrac{1}{h_{r}^{\beta }}\sum_{k\in I_{r}}\left[ f\left( d\left( \Delta
^{m}X_{k},X_{0}\right) \right) \right] ^{p_{k}} &=&\dfrac{1}{h_{r}^{\beta }}%
\sum_{1}\left[ f\left( d\left( \Delta ^{m}X_{k},X_{0}\right) \right) \right]
^{p_{k}}+\dfrac{1}{h_{r}^{\beta }}\sum_{2}\left[ f\left( d\left( \Delta
^{m}X_{k},X_{0}\right) \right) \right] ^{p_{k}} \\
&\leq &\dfrac{1}{h_{r}^{\beta }}\sum_{1}\max \left( K^{h},K^{H}\right) +%
\dfrac{1}{h_{r}^{\beta }}\sum_{2}\left[ f\left( \varepsilon \right) \right]
^{p_{k}} \\
&\leq &\max \left( K^{h},K^{H}\right) \dfrac{1}{h_{r}^{\beta }}\left\vert
\left\{ k\in I_{r}:f\left( d\left( \Delta ^{m}X_{k},X_{0}\right) \right)
\geq \varepsilon \right\} \right\vert \\
&&+\dfrac{h_{r}}{h_{r}^{\beta }}\max \left( f\left( \varepsilon \right) ^{h},%
\text{ }f\left( \varepsilon \right) ^{H}\right) .
\end{eqnarray*}%
Hence we get $X\in w_{p}^{\beta }\left( \theta ,f,F,\Delta ^{m}\right) $ for 
$r\in 
\mathbb{N}
$ from above inequalities.

\bigskip

\noindent \textbf{Theorem 3.4. }If $\lim p_{k}>0$ and sequence $X=\left(
X_{k}\right) $ is strongly $w_{p}^{\beta }\left( \theta ,f,F,\Delta
^{m}\right) -$summable to $X_{0},$ then $w_{p}^{\beta }\left( \theta
,f,F,\Delta ^{m}\right) -\lim X_{k}=X_{0}$ is unique.

\textbf{Proof. }Let\textbf{\ }$\lim p_{k}=s>0,$ $w_{p}^{\beta }\left( \theta
,f,F,\Delta ^{m}\right) -\lim X_{k}=X^{\prime }$ and $w_{p}^{\beta }\left(
\theta ,f,F,\Delta ^{m}\right) -\lim X_{k}=X^{\prime \prime }$. Then we
obtain

\begin{equation*}
\lim_{r\rightarrow \infty }\dfrac{1}{h_{r}^{\beta }}\sum_{k\in I_{r}}\left[
f\left( d\left( \Delta ^{m}X_{k},X^{\prime }\right) \right) \right]
^{p_{k}}=0
\end{equation*}%
and

\begin{equation*}
\lim_{r\rightarrow \infty }\dfrac{1}{h_{r}^{\beta }}\sum_{k\in I_{r}}\left[
f\left( d\left( \Delta ^{m}X_{k},X^{\prime \prime }\right) \right) \right]
^{p_{k}}=0
\end{equation*}%
On the other hand, since $f$ is increasing and $f\left( x+y\right) \leq
f\left( x\right) +f\left( y\right) $ we can write

\begin{eqnarray*}
\dfrac{1}{h_{r}^{\beta }}\sum_{k\in I_{r}}\left[ f\left( d\left( X^{\prime
},X^{\prime \prime }\right) \right) \right] ^{p_{k}} &\leq &\dfrac{D}{%
h_{r}^{\beta }}\sum_{k\in I_{r}}\left( \left[ f\left( d\left( \Delta
^{m}X_{k},X^{\prime }\right) \right) \right] ^{p_{k}}+\left[ f\left( d\left(
\Delta ^{m}X_{k},X^{\prime \prime }\right) \right) \right] ^{p_{k}}\right) \\
&\leq &\dfrac{D}{h_{r}^{\beta }}\sum_{k\in I_{r}}\left[ f\left( d\left(
\Delta ^{m}X_{k},X^{\prime }\right) \right) \right] ^{p_{k}}+\dfrac{D}{%
h_{r}^{\beta }}\sum_{k\in I_{r}}\left[ f\left( d\left( \Delta
^{m}X_{k},X^{\prime \prime }\right) \right) \right] ^{p_{k}}\rightarrow 0,
\end{eqnarray*}%
where $p_{k}=H$ and $D=\left( 1,2^{H-1}\right) $. Left right of above
inequalities tends to zero as $r\rightarrow \infty $ and because of $\lim
p_{k}=s$ we obtain%
\begin{equation*}
\lim_{r\rightarrow \infty }\left[ f\left( d\left( X^{\prime },X^{\prime
\prime }\right) \right) \right] ^{p_{k}}=\left[ f\left( d\left( X^{\prime
},X^{\prime \prime }\right) \right) \right] ^{s}
\end{equation*}%
Hence $X^{\prime }=X^{\prime \prime }$ and the limit is unique.

\end{document}